\documentclass{amsart}
\usepackage{amssymb,amscd,amsthm,verbatim}
\usepackage{graphicx}

\newtheorem{proposition}{Proposition}

\newtheorem{theorem}[proposition]{Theorem}
\newtheorem{prop}[proposition]{Proposition}
\newtheorem{corollary}[proposition]{Corollary}
\newtheorem{lemma}[proposition]{Lemma}
\newtheorem{conjecture}[proposition]{Conjecture}

\newtheorem{remark}[proposition]{Remark}

\theoremstyle{definition}
\newtheorem{defn}[proposition]{Definition}

\def\iin{{\mathrm{in}_\preceq}}
\def\Hilb{{\mathrm{Hilb}}}
\def\sd{{\text{sd}}}
 \def\conv{{\text{conv}}}

\def\ZZ{{\mathbb{Z}}}
\def\RR{{\mathbb{R}}}
\def\NN{{\mathbb{N}}}
\def\CC{{\mathbb{C}}}
\def\KK{{\mathbb{K}}}
\def\kk{{{\Theta}}}

\def\CR{{\mathcal{R}}}

\newcommand{\link}{\mathrm{link}}

\newcommand{\F}{\mathcal F}

\newcommand{\lt}{\mathrm{lm}_\preceq}

\numberwithin{equation}{section}

\newcommand{\D}{\mathcal D}

 \DeclareMathOperator{\mot}{mod}

\begin{document}

\title[Type-B generalized triangulations and 
determinantal ideals]
{Type-B generalized triangulations and 
determinantal ideals}

\author[Daniel Soll]{Daniel Soll}
\email{daniel.soll@googlemail.com}

\author[Volkmar Welker]{Volkmar Welker}
\email{welker@mathematik.uni-marburg.de}

\address{Philipps-Universit\"at Marburg\\
         Fachbereich Mathematik und Informatik\\
         35032 Marburg. Germany}


\thanks{Both authors were supported by EU Research Training Network
``Algebraic Combinatorics in Europe'', grant HPRN-CT-2001-00272}

\begin{abstract}
For $n\geq 3$, let
$\Omega_n$ be the set of line segments between the vertices of a convex $n$-gon.
For $j\geq 2$, a $j$-crossing is a set of $j$ line segments pairwise
intersecting in the relative interior of the $n$-gon. For $k\geq 1$, let
$\Delta_{n,k}$ be the simplicial complex of (type-A) generalized triangulations, i.e.
the simplicial complex of subsets of $\Omega_n$ not containing
any $(k+1)$-crossing.
 
The complex $\Delta_{n,k}$ has been the central object of numerous papers. 
Here we continue this work by considering the complex of type-B generalized
triangulations. For this we identify line-segments in $\Omega_{2n}$ which
can be transformed into each other by a $180^\circ$-rotation of the
$2n$-gon. Let $\F_n$ be the set $\Omega_{2n}$ after identification, then the
complex $\D_{n,k}$ of type-B generalized triangulations
is the simplicial complex of subsets of $\F_n$ not containing
any $(k+1)$-crossing in the above sense. For $k = 1$, we have that $\D_{n,1}$ is
the simplicial complex of type-B triangulations of the $2n$-gon as defined in
\cite{Si} and decomposes into a join of an 
$(n-1)$-simplex and the boundary of the $n$-dimensional cyclohedron. We
demonstrate that $\D_{n,k}$ is a pure, $k(n-k)-1+kn$ dimensional complex
that decomposes into a $kn-1$-simplex and a $k(n-k)-1$ dimensional homology sphere.
For $k=n-2$ we show that this homology-sphere is in fact the boundary of a
cyclic polytope. We provide a lower and an upper bound for the number of
maximal faces of $\D_{n,k}$.

On the algebraical side we give a term-order on the monomials in the variables $X_{ij}, 1\leq i,j\leq n$,
such that the corresponding initial ideal of the determinantal ideal
generated by the $(k+1)$ times $(k+1)$ minors of the generic $n \times n$
matrix contains the Stanley-Reisner ideal of $\D_{n,k}$. We show that the minors
form a Gr\"obner-Basis whenever $k\in\{1,n-2,n-1\}$ thereby proving the
equality of both ideals and the
unimodality of the $h$-vector of the determinantal ideal in these cases.
We conjecture this result to be true for all values of $k<n$.

\end{abstract}

\maketitle

\section{Introduction and Basic Definitions}

Generalized associahedra have been subject to fruitful and
intensive study recently (see for example \cite{FZ}). In this
research associahedra are defined uniformly for all root systems --
the classical associahedron being the type-A case. In a second
stream originating in the work of Dress, Koolen \& Moulton \cite{DGJM}
and Nakamigawa \cite{Na} (see also \cite{DKKM}, \cite{Jo}, \cite{JW}) a
generalization of the type-A associahedron into a different 
direction has been shown to exhibit very nice 
combinatorial, geometric and algebraic properties. 

In this paper we carry out an analogous generalization for the type-B 
associahedron -- the cyclohedron (see for example \cite{Si}, \cite{Re}).
We verify similar nice combinatorial, geometric and algebraic properties.
In particular, we show that it is a homology sphere, give bounds on the number
of facets and relate it to determinantal ideals.

Note, that by polytope duality there are two polytopes -- one simple, one
simplicial -- that can be called associahedron (resp. cyclohedron). In this paper, by
associahedron (resp. cyclohedron) we will always mean the simplicial polytope.

The paper is organized as follows. This section gives a brief introduction to the
topic and provides basic definitions on simplicial complexes
and Gr\"obner bases.

In Section \ref{typea} we recall the results on the generalized type-A 
associahedron. 
In Section \ref{typeb} we then state the main results of this paper. 
In Sections \ref{proofdiameters}, \ref{proofhomsphere} and \ref{definition} we then provide the 
proofs of the main results. 

In order to recall known facts and formulate our own results 
we need to introduce basic notions
about simplicial complexes and Gr\"obner bases.

A simplicial complex on vertex-set $V$ is a set-system $\Delta\subset 2^V$
such that $\sigma\in \Delta$ implies $\tau\in\Delta$ for all
$\tau\subset\sigma$. Elements of $\Delta$ are called faces. The dimension 
$\dim \sigma$ of
a face $\sigma\in\Delta$ is the number of elements in $\sigma$ reduced by
one. 
The dimension  $\dim \Delta$ of the complex $\Delta$ is the maximal dimension of a facet in
$\Delta$. If all facets have the same dimension the complex is said to
be pure. The $f$-vector of a $(d-1)$-dimensional complex $f(\Delta)=(f_{-1},f_0,\dots,f_{d-1})$
counts the faces by their dimension, i.e. $f_i$ is the number of faces
of dimension $i$. The $h$-vector is a transformation
of the $f$-vector defined as follows.
We set $h(t) := f(t-1) = \sum_{i=0}^dh_it^{d-i}$ -- the $h$-polynomial -- for 
$f(t):=\sum_{i=0}^{d-1}f_{i-1}t^{d-i}$ -- the $f$-polynomial -- and obtain the $h$-vector by
$h(\Delta) = (h_0, \ldots, h_d)$. 
A simplicial complex may be defined by its
inclusion-maximal faces, the facets, or by its inclusion-minimal nonfaces.
The ideal $I_\Delta$ in $T = k[x_v~|~v \in V]$ generated by $\prod_{v \in sigma} x_v$
for the inclusion minimal non-faces $\sigma$ is called Stanley-Reisner ideal of
$\Delta$ and $k[\Delta] := T/I_\Delta$ the Stanley-Reisner ring of $\Delta$. 
The Krull-dimension of $k[\Delta]$ is well known to equal $\dim \Delta + 1$ and 
the entry $f_{d-1}$ for $d-1 = \dim \Delta$ is known to be the multiplicity of 
$k[\Delta]$ (see \cite[Chapter 5]{BH} for more details) . 
More generally, the Hilbert-series of $k[\Delta]$ for a $d$-dimensional 
simplicial complex $\Delta$ with $h$-vector $(h_0, \ldots, h_d)$ is given by
$$\Hilb(k[\Delta],t) = \frac{h_0 + \cdots + h_dt^d}{(1-t)^{d+1}}.$$ 
For simplicial complexes $\Delta,\Sigma$ over disjoint vertex-sets $V_1,V_2$ one
defines the join $\Delta\star\Sigma$ as the simplicial complex on vertex-set $V_1\cup V_2$
and simplices $\tau=\delta\cup\sigma$ for
$\delta\in\Delta,\sigma\in \Sigma$.
The minimal non-faces of $\Delta * \Sigma$ are the minimal non-faces of $\Delta$
and of $\Sigma$. In particular, $k[\Delta * \Sigma] = 
k[\Delta] \otimes k[\Sigma]$. Note, that in the special case when $\Sigma$ is
simplex all minimal non-faces of $\Delta * \Sigma$ are minimal non-faces of
$\Delta$. This case will become crucial later in this work. 

Since the Stanley-Reisner ideals are monomial ideals 
(i.e., ideals generated by monomials), we will need Gr\"obner basis theory 
and initial ideals for establishing the link
between the Stanley-Reisner ideals of our complexes and determinantal ideals 
in Section \ref{typeb}. 
In the following paragraph we introduce the basic terminology of
Gr\"obner basis theory in our setting (for more details see \cite{AL}).

Let $\preceq$ be a term order on the monomials in the polynomial ring $S = k[x_1, \ldots, x_n]$. 
For a polynomial $f\in S$ we let
$\lt(f)$ be the leading monomial of $f$, i.e. the largest monomial appearing with
non-zero coefficient in $f$ with respect to $\preceq$.
For a subset $M\subset S$ let $\iin (M)$ be the ideal generated by
the set of monomials $\{\lt(f)\big| f\in M\}$. Recall that a Gr\"obner-basis $G$ of an ideal
$I$ is a subset $G \subseteq I$ such that $\iin(G)=\iin(I)$.
It is well known that $\iin(I)$ and $I$ share their Hilbert-series. In particular,
their Krull-dimension and multiplicity coincide (see \cite{AL}).

\section{Type-A $k$-Triangulations} \label{typea}

The definition of type-A $k$-triangulations takes advantage of geometric
properties of the convex $n$-gon, $n\geq 3$, which we define as the convex hull 
of the $n$-th roots of unity $\{1=\xi_0,\dots,\xi_{n-1}\}$ in $\CC \cong \RR^2$
numbered in clockwise order.
A diagonal between the $i$th and $j$th root of unity is the line segment 
$\partial_{i,j}:=\{\lambda\xi_i+(1-\lambda)\xi_j\big|\lambda\in
(0,1)\}\subset \RR^2$ which we may identify with the set $\{i,j\}$ and
often denote as $ij$. Let $\Omega_{n}=\{ij\subset [n],\,i\not=j\}$ be the
set of diagonals of the $n$-gon. A $(k+1)$-crossing is 
a $(k+1)$-set of diagonals mutually intersecting in the 
relative interior of the $n$-gon. As a consequence, 
diagonals in the set $\Gamma_{n,k}:=\{ij\big|n-k\leq |j-i|\leq k\}$ cannot be part of a
$(k+1)$-crossing. For $n\geq 2k+1$ and $k\geq 1$ one defines
$\Delta_{n,k}$ as the
simplicial complex on vertex-set $\Omega_n$ and the minimal nonfaces being
the set of all $(k+1)$-crossings. When $k=1$, the facets of this complex correspond to
triangulations of the $n$-gon. It is easy to see that $\Delta_{n,k}$ can be
decomposed as a join $\Delta_{n,k}=:\Delta^\ast_{n,k}\star 2^{\Gamma_{n,k}}$,
where the complex $\Delta^\ast_{n,k}$ is set of all faces of $\Delta_{n,k}$ contained in 
$\Omega_{n,k} := \Omega_n \setminus \Gamma_{n,k}$. 
Following \cite{Jo}, we call facets of both complexes $\Delta_{n,k}$ and
$\Delta_{n,k}^\ast$ (type-A)-\textsc{generalized
$k$-triangulations}. The prefix type-A is motivated by the theory of
cluster-complexes where the complex $\Delta_{n,1}$ is shown to be the
cluster-complex of an arbitrary cluster-algebra of type $A_n$, while the
complex $\D_{n,1}$ to be defined in the next section is the cluster-complex
of type $B_n$ cluster-algebras. For more
information we refer the reader to \cite{FR} and \cite{FZ}.

Both complexes $\Delta_{n,k}$ and $\Delta_{n,k}^\ast$
have been the central objects of several papers. In the sequel we will list
some more recent results, most of them originating in work of Dress, Koolen and Moulton \cite{DKM1}, and refer the reader to \cite{CP} for prior developments.

\begin{figure}
\includegraphics[scale=0.5]{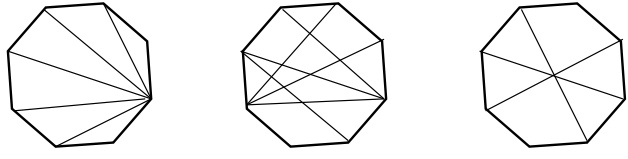}\
\caption{Generalized type-A $k$-triangulations for $n=8, k=1,2,3$.}\label{ufacettenn8k123}
\end{figure}

\begin{theorem}[\cite{DKM1},\cite{Na}] \label{dim} The complex
$\Delta_{n,k}^\ast$ is pure of dimension
$k(n-2k-1)-1.$\end{theorem}

\begin{theorem}[\cite{Jo}]\label{anzfac}
The number of facets of $\Delta_{n,k}$ resp. $\Delta_{n,k}^\ast$ is given by the following determinant:
\begin{eqnarray*}
\det\left(\begin{array}{ccccc}C_{n-2}&
C_{n-3}&\dots&C_{n-k}&C_{n-k-1}\\
C_{n-3}& C_{n-4}&\dots&C_{n-k-1}&C_{n-k-2}\\
\dots&\dots&\dots&\dots&\dots\\
C_{n-k-1}& C_{n-k-2}&\dots&C_{n-2k+1}&C_{n-2k}
\end{array}\right)&=&\prod_{1\leq i\leq j\leq n-2k-1}\frac{i+j+2k}{i+j}
\end{eqnarray*}
where $C_n=\frac{1}{n+1}{2n\choose n}$ is the $n$th Catalan-Number.
\end{theorem}

Recently Krattenthaler \cite{Kr} gave an alternative
proof of Theorem \ref{anzfac} in the context of Fomin's growth diagrams
(see also \cite{dM} for related work).

\begin{theorem}[\cite{DGJM}]\label{sphere} The geometric realization of the 
simplicial complex $\Delta_{n,k}^\ast$ is a piecewise linear sphere.
\end{theorem}

Since every boundary complex of a simplicial polytope is a piecewise linear
sphere, Theorem \ref{sphere} has raised the question whether $\Delta^\ast_{n,k}$ is
polytopal, which is commonly
believed to have a positive answer. 

The question of polytopality takes its appeal  not only from
geometry, but also from enumerative combinatorics. An affirmative answer would reveal 
certain properties of the
$h$-vector of $\Delta_{n,k}$ as they are stated in the famous $g$-theorem by
Billera, Lee, McMullen and Stanley, see \cite{Hib}. A consequence of special
interest in enumerative combinatorics is the unimodality ($h_0\leq h_2\leq
h_{\lfloor\frac{d}{2}\rfloor}=h_{\lfloor\frac{d+1}{2}\rfloor}\geq\dots h_d)$
of the $h$-vector.

In special cases the polytopality has been ascertained:

\begin{prop}\cite{DGJM}\label{polytopaltypa} The complex $\Delta_{n,k}^\ast$ is the boundary-complex of a simplicial
polytope for $2k+1\leq n\leq 2k+3$ and for $k=1$, i.e. it . If $k=1$, this polytope is the well known associahedron, for
$n=2k+1$ it is the $(-1)$-sphere, for $n=2k+2$ it is the $k$-simplex and for
$n=2k+3$ it is a cyclic polytope.
\end{prop}

Recall that the cyclic polytope $\mathcal C_d(n)$ is the convex hull of $n$ different points on the moment curve
$M_d:=\{(1,t,t^2,\dots,t^d),\,t\in \RR\}$. It is well known that the combinatorics of its boundary-complex 
does not depend on the choice of the points.
Polytopal realizations of the associahedron can be found in \cite{Ha}, \cite{Le} and \cite{Lo}.

Finally, there is an unexpected relation of $\Delta_{n,k}$ and the ideal
$P_{n,k}$ of Pfaffians of degree $k+1$ of a generic $n \times n$-skew 
symmetric matrix in the indeterminates $x_{ij}$, $1 \leq i < j \leq n$. 
For that we identify the variable in the Stanley-Reisner ring of
$\Delta_{n,k}$ indexed with the diagonal $\{i,j\}$ with the variable
$x_{ij}$.

\begin{theorem}[\cite{JW}] \label{pfaffian} 
For $2 \leq 2k + 2 \leq n$ there is a term order $\preceq$ for which 
$\iin (P_{n,k}) = I_{\Delta_{n,k}}$.
\end{theorem} 
 
\section{Main results and type-B $k$-triangulations} \label{typeb}

For the type-B case we identify line-segments in $\Omega_{2n}$ which
can be transformed into each other by a $180^\circ$-rotation of the
$2n$-gon. We will write $\bar M$ for a rotated set
$M\subset\Omega_{2n}$ and $\bar d$ for the rotated $d\in\Omega_{2n}$.

Line-segments with $d=\bar d$, i.e. those that cross the origin in the
$2n$-gon, will be called diameters.
Let $\F_n$ be the set $\Omega_{2n}$ after identification, then for $1\leq k\leq
n-1$ the complex of type-B generalized $k$-triangulations
$\D_{n,k}$ is the simplicial complex of subsets of $\F_n$ not containing
any $(k+1)$-crossing in the above sense.

When $k=1$, the facets of this complex correspond to type-B
triangulations of the $2n$-gon as defined in \cite{Si}.

\begin{figure}
\includegraphics[scale=0.5]{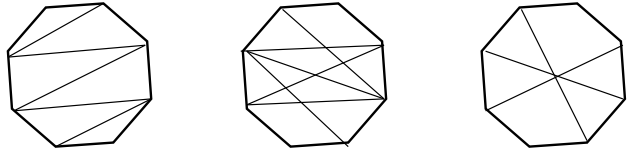}\
\caption{Generalized type-B $k$-triangulations for $n=4, k=1,2,3$.}\label{facettenn8k123}
\end{figure}

Let $\F_{n,k}$ be the set of classes in $\Omega_{2n,k}$ and $\kk_{n,k}$
the set of classes in $\Gamma_{2n,k}$ with respect to the identification
mentioned above. Again it is easy to see that $\D_{n,k}$ can be
written as a join $\D_{n,k}=:\D^\ast_{n,k}\star 2^{\kk_{n,k}}$,
where the complex $\D^\ast_{n,k}$ is defined on the vertex-set
$\F_{n,k}$. Facets of $\D_{n,k}^\ast$ will be called type-B generalized $k$-triangulations.
See Figure \ref{facettenn8k123} for some examples.

We will also consider the complex $\Delta^{\ast\text{symm}}_{2n,k}$ that
we define to be the subcomplex of $\Delta^\ast_{2n,k}$ which is generated by the
facets $F$ of $\Delta_{2n,k}$ with $F=\bar F$.

We are now in the position to state our main results.

\begin{theorem}\label{diameters}
For all $1\leq k\leq n-1$  each facet of $\Delta^{\ast\text{symm}}_{2n,k}$ contains
exactly $k$ diameters.
\end{theorem}

As an instant consequence of the preceding theorem and Theorem \ref{dim} we get:

\begin{corollary}\label{dimtypb}
The complex $\D^\ast_{n,k}$ is pure of dimension $k(n-k)-1$.
\end{corollary}

Concerning polytopality we prove a type-B analogue of Proposition 
\ref{polytopaltypa}:

\begin{prop}\label{polytopaltypb} The complex $\D_{n,k}^\ast$ is the boundary-complex of a simplicial
polytope for $k\in\{1,n-2,n-1\}$. If $k=1$, this polytope is the 
cyclohedron (\cite{Si}), for
$k=n-1$ it is the $k$-simplex and for $k=n-2$ it is a cyclic polytope.
\end{prop}

The cyclohedron was first introduced by \cite{St}. For an explicit realization as a polytope see \cite{HL}.
This plies us with the formulation of the following conjecture.

\begin{conjecture}\label{polytopality}
The complex $\D^\ast_{n,k}$ is polytopal for all $1\leq k\leq n-1$.
\end{conjecture}

We cannot prove Conjecture \ref{polytopality} in its full strength.
But as a step towards the conjecture we prove the following theorem.

\begin{theorem}\label{homsphere} For $1\leq k\leq n-1$ the complex 
$\D^\ast_{n,k}$ is a $\mot 2$-homology-sphere.
\end{theorem}

For a prime $p$, we call a simplicial complex $\Delta$ a $\mot p$-homology-sphere if
for each face $\sigma\in \Delta$ we have 
$$H_\star(\link_{\Delta}(\sigma),\ZZ_2)\cong H_\star(S^{k(n-k)-1-\dim \sigma -1},\ZZ_2).$$
In particular, for $\sigma = \emptyset$ we then have 
$$H_\star(\link_\Delta \emptyset ,\ZZ_2) = H_\star(\Delta ,\ZZ_2)\cong H_\star(S^{k(n-k)-1},\ZZ_2)$$.

Theorem \ref{homsphere} provides us with the prerequisites in order to obtain
bounds on the number of generalized type-B triangulations.
While the upper bound is an immediate consequence of \cite{No1},\cite{No2}
the proof of the lower bound can be found in Section \ref{definition}.

\begin{theorem}\label{bounds}
For all $1\leq k\leq n-1$ the number $T(n,k)$ of type-B generalized
$k$-Triangulations of the $2n$-gon satisfies
\begin{itemize}

\item  the lower bound

\begin{eqnarray}
T(n,k)&\geq&\det\left[{2n-i-j \choose n-i}\right]_{i,j=1,\dots,k}\label{ht1}\\
&=&\det\left[{2(n-k) \choose n-k+i-j}\right]_{i,j=1,\dots,k}\nonumber\\
&=& \prod_{h=1}^{n-k}\prod_{i=1}^{k}\prod_{j=1}^{n-k}\frac{h+i+j-1}{h+i+j-2}\nonumber 
\end{eqnarray}

\item the upper bound given by the number of facets of the $k(n-k)$ dimensional
cyclic polytope with $n(n-k)$ vertices
$$T(n,k)\leq\begin{cases}2\cdot{(n-k)^2+\frac{k(n-k)-1}{2}\choose
(n-k)^2}&\text{if }k(n-k)\textrm{ odd,}\\
2\cdot{(n-k)^2+\frac{k(n-k)}{2}\choose
(n-k)^2}+{(n-k)^2-1+\frac{k(n-k)}{2}\choose(n-k)^2-1}&\text{if }
k(n-k)\textrm{ even.}\end{cases}$$

\end{itemize}
\end{theorem}

\begin{remark}
If we fix $l=(n-k)$ and let $n$ and therefore $k$ go to infinity, the quotient
of the lower and the upper bound converges to $\frac{2!\cdot 3!\cdots
(l-1)!\cdot l^2!}{2\cdot l!\cdot(l+1)!\cdots (2l-1)!\cdot(l/2)^{l^2}}$.
\end{remark}

It has been shown in \cite{Si} that $T(n,1)={2n-2\choose n-1}$, which means
that the lower bound is met in this case, which is also true when
$\D_{n,n-1}$ is the simplex of dimension $(n-1)$. Proposition \ref{polytopaltypb}
yields that the lower bound and the upper bound coincide in the
case $k=n-2$. This and results in Section \ref{definition} lead us to the next conjecture.

\begin{conjecture}\label{Vermutung1}
The inequality (\ref{ht1}) is an equality. 
\end{conjecture}

Finally, we provide a description of the Stanley-Reisner ideal of $\D_{n,k}$ in terms of submatrices. This will lead to a partially conjectural connection 
of $\D_{n,k}$ and the ideal of minors of degree $k+1$ of a generic $n \times
n$ matrix, analogous to Theorem \ref{pfaffian}. 
In this case the identification of the set of variables used in the 
Stanley-Reisner ideal of $\D_{n,k}$ and the set of variables of a generic
$n \times n$ matrix is not as obvious as in the Pfaffian case. 
Therefore, for the description of $I_{D_{n,k}}$ in terms of the entries of a
generic $n \times n$ matrix, we apply the following 
bijections between the ground set of 
$D_{n,k}$ and the index pairs $(i,j)$ for $1 \leq i,j \leq n$. 

\begin{eqnarray*}
\Psi:\quad\{\{a_1<b_1\},\{a_2<b_2\}\}\mapsto \begin{cases}((a_i+1)\mot n,b_i\mot
n)&\textrm{ if }b_i-a_i\leq n\\((b_i+1)\mot n,a_i\mot n)&\textrm{ if
}b_i-a_i>n,\end{cases}
\end{eqnarray*}

The map $\Psi$ will be studied in more detail in the proof of Theorem \ref{Stan}.
Let $1\leq k < n$.
For nonempty subsets 
$A,B\subset [n]$ with $\sharp A=\sharp B=(k+1)$ we
define $M(A,B):=(x_{ij})_{i\in A,j\in B}$ to be the corresponding row- and column-selected submatrix.
In the following we always assume $k\mot k=k$. We let
\begin{eqnarray}
N(A,B)&:=&\{\left(a_{(i+l)\mot (k+1)},b_i\right),i=1,\dots,(k+1)\}\label{nichtseiten}
\end{eqnarray}
where $(k+1)\geq \ell$ is chosen such that the following two conditions hold:
\begin{eqnarray}
a_{i+\ell}&>&b_i\textrm{ for all }i=1,\dots,(k+1)-l \\
a_{j+\ell-1}&\leq&b_j\textrm{ for one }j\in\{1,\dots,(k+1)-l+1\}\textrm{ or }l=0.\end{eqnarray}
For $k+1=1$ we let $N(A,B)=\{(a_1,b_1)\}$.
(See Figure \ref{submatrix} for an example.) The set $N(A,B)$ consists of all
indices of matrix entries on the longest diagonal of $M(A,B)$ that lies strictly below the
main-diagonal of $X$ augmented by the indices of the entries on the 
complementary diagonal.

\setlength{\unitlength}{0.75mm}

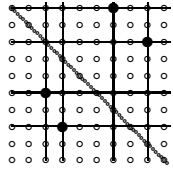
\begin{figure}\label{submatrix}
\begin{center}
\begin{picture}(32,32)(-2,-2)

\multiput(0,0)(3,0){10}{\circle{1}}
\multiput(0,3)(3,0){10}{\circle{1}}
\multiput(0,6)(3,0){10}{\circle{1}}
\multiput(0,9)(3,0){10}{\circle{1}}
\multiput(0,12)(3,0){10}{\circle{1}}
\multiput(0,15)(3,0){10}{\circle{1}}
\multiput(0,18)(3,0){10}{\circle{1}}
\multiput(0,21)(3,0){10}{\circle{1}}
\multiput(0,24)(3,0){10}{\circle{1}}
\multiput(0,27)(3,0){10}{\circle{1}}

\put(6,0){\line(0,1){28}}
\put(9,0){\line(0,1){28}}
\put(18,0){\line(0,1){28}}
\put(24,0){\line(0,1){28}}

\put(0,6){\line(1,0){28}}
\put(0,12){\line(1,0){28}}
\put(0,21){\line(1,0){28}}
\put(0,27){\line(1,0){28}}

\put(9,6){\circle*{2}}
\put(6,12){\circle*{2}}
\put(24,21){\circle*{2}}
\put(18,27){\circle*{2}}

\multiput(0,27)(.5,-.5){56}{\circle{0.2}}

\end{picture}
\caption{$N(A,B)$ for
$n=10,\,k=3,\,A=\{1,3,6,8\},\,B=\{3,4,7,9\}$}
\end{center}

\end{figure}

We then have

\begin{theorem}\label{Stan}
The Stanley-Reisner ideal of $\D^\ast_{n,k}$
is generated by all monomials
$$m_{A,B}=\prod_{(i,j)\in N(A,B)}x_{ij}\text{~such that~} \sharp A=\sharp
B=(k+1).$$
\end{theorem}

The preceding result allows us to connect the complex $\D^\ast_{n,k}$ to 
the determinantal ideal $I_{n,k}$. 
For $1 \leq k<n$ the ideal $I_{n,k}$ is defined as follows. 
Let $S:=\KK[x_{11},\dots,x_{nn}]$ be the polynomial ring over the field $\KK$
and consider the $x_{ij}$ as entries in the $(n \times n)$-matrix $X$ of indeterminates. 
Then $I_{n,k}$ is the ideal
generated by the $(k+1)$-minors, i.e. by the formal determinants of all square-submatrices
of size $(k+1)\times (k+1)$ of $X$.
The determinantal ideal $I_{n,k}$ is one of the classical objects in commutative algebra
which offers numerous links to the theory of invariants, group representation
theory and combinatorics. For more details see \cite{BV} and \cite{BC}.

In Section \ref{definition} we will establish a term order $\preceq$ for which 
$I_{D_{n,k}} \subseteq \iin (I_{n,k})$ (see Theorem \ref{NAB}). 
Indeed, we able to verify that for $I = I_{D_{n,k}}$ and $J = \iin (I_{n,k})$
all assumptions of the following lemma are satisfied except for (ii).
The lemma was first successfully applied in \cite{JW}
for the proof of Theorem \ref{pfaffian}. A much more general version
can be found in \cite[Lemma 4.2]{CHT}.

\begin{lemma}[\cite{JW}] \label{monomialinclusion}
 Let $T = k[x_1, \ldots, x_\ell]$ be the polynomial ring in $\ell$
 variables. Suppose that $I \subseteq J$ are monomials ideals in $T$
 such that the following hold:
   \begin{itemize}
     \item[(i)] $\dim(T/I) = \dim(T/J)$.
     \item[(ii)] $e(T/I) = e(T/J)$.
     \item[(iii)] $I = I_\Delta$ for a pure simplicial complex $\Delta$
                  on ground set $[\ell]$.
   \end{itemize}

   Then $I = J$.
\end{lemma}

Thus in order to establish equality between $I = I_{D_{n,k}}$ and $J= \iin (I_{n,k})$
it remains to verify condition (ii). Even though we are not able to prove this fact, we can
deduce one inequality from the fact that $I \subseteq J$. 
 
\begin{lemma}\label{halb}
Let $\Delta$ be a pure simplicial complex over the vertex-set $\{x_1,\dots,x_n\}$
and let $J$ be a monomial ideal in $S=\KK[x_1,\dots,x_n]$.
If the Stanley-Reisner ideal $I_\Delta$ is contained in $J$ and their
Krull-dimensions coincide, then for the multiplicities holds
$$e(S\slash I_\Delta)\geq e(S\slash J).$$
\end{lemma}

\begin{proof}{}
We consider the polarization $J^{pol}$ of $J$, i.e. we replace every occurrence
of $x_i^{\alpha_i}$ in a generator of $J$ by additional variables $x_i^1,\dots,x_i^{\alpha_i}$.
Let $J^{pol}$ be the ideal generated by the replaced variable in the
ring $S'$, which we may write as $k[x_1,\dots,x_{n+m}]$ after a suitable renaming
of the variables.
The ideal $J^{pol}$ is squarefree, therefore exists a simplicial complex
$\Gamma$, such that for the Stanley-Reisner ideal
$I_\Gamma$ holds: $J^{pol}=I_\Gamma$.
Now we take the generation of
$I_\Delta$ in $S'$; it corresponds to a squarefree ideal $I'$, that can be
interpreted as the Stanley-Reisner ideal of the complex $\Delta\star
2^{\{x_{n+1},\dots,x_{n+k}\}}$. (Where we write
$2^{\{x_{n+1},\dots,x_{n+k}\}}$ for the simplex with vertex-set
$\{x_{n+1},\dots,x_{n+k}\}$). All in all we obtain for both
Stanley-Reisner ideals:
$$I_{\Delta\star 2^{\{x_{n+1},\dots,x_{n+k}\}}}\subset I_\Gamma.$$
We can conclude for the corresponding simplicial complexes:
$$\Gamma\subset\Delta\star2^{\{x_{n+1},\dots,x_{n+k}\}}\textrm{ and
}\dim\Gamma=\dim \Delta\star2^{\{x_{n+1},\dots,x_{n+k}\}}.$$
Therefore the number of faces of maximal dimension of $\Gamma$ is smaller or equal than
the number of facets of
$\Delta\star2^{\{x_{n+1},\dots,x_{n+k}\}}$, and therefore smaller or equal
to the number of facets of $\Delta$. Since $\Delta$ is a pure complex, this
number is the multiplicity of $\Delta$ which equals the multiplicity of
$S\slash I_{\Delta}$.
\end{proof}

For the special cases  $k\in\{1,n-2,n-1\}$ we indeed can establish 
assumption (ii) of Lemma and hence show that 
$I_{D_{n,k}} = \iin (I_{n,k})$. 
Moreover, the term-order $\preceq$ satisfies 
the following conditions:

\begin{itemize}
 \item[(S)] The initial ideal $\iin (I_{n,k})$ is the
            Stanley-Reisner ideal of the simplicial complex $\D_{n,k}$ which
            decomposes into a join $\D^\ast_{n,k} \star 2^\kk$ of a
            triangulation $\D^\ast_{n,k}$ of a sphere and a full simplex $2^\kk$
            on ground set $\kk$.
 \item[(M)] The cardinality of the set $\kk$ is the absolute value
            of the $a$-invariant of the quotient of the polynomial
            ring by $I_{n,k}$.
\end{itemize}

Term-orders and ideals that satisfy these two condition have recently appeared
in several places in the literature (see \cite{CHT} for some original results in
this direction and an exhaustive survey of the previous known instances).
For the type-A case of our situation this are the results from \cite{JW}.
Initial ideals which satisfy (S) and (M) 
will be called spherical initial ideals. 

Thus the a verification of Conjectures \ref{polytopality} and 
\ref{Vermutung1} would imply: 

\begin{conjecture} 
For the term-order $\preceq$ defined in Section \ref{definition} 
the initial ideal of $I_{n,k}$ is a spherical initial ideal. 
\end{conjecture}

With the help of Lemma \ref{monomialinclusion} we get the following equivalent 
formulations of Conjecture \ref{Vermutung1}.

\begin{lemma} \label{dreitenoere}
The following statements are equivalent:
\begin{itemize}
\item The $(k+1)$-minors form a Gr\"obner-bases for the term-order $\preceq$.
\item The Stanley-Reisner Ideal of $\D_{n,k}$ coincides with $\iin(I_{n,k})$.
\item The number of generalized type-B triangulations is counted according
to Conjecture \ref{Vermutung1}.
\end{itemize}
\end{lemma}

\section{Proof of Theorem \ref{diameters}} \label{proofdiameters}

To determine the dimension and ascertain the pureness of $\D_{n,k}$ we
symmetrize an approach of Jonsson \cite{Jo}. The original formulation in
terms of triangulations was present in a preprint version \cite{Jopre} 
of \cite{Jo} but 
was replaced by a formulation in terms of polyominos.
Indeed, all the results mentioned below are special cases of results 
from \cite{Jo}, but explicitly stated only in the preprint \cite{Jopre}. 

In this section we will often consider the complexes $\Delta_{n,k}$ and 
$\Delta_{n-1,k}$ in parallel. It turns out that in the proofs it is 
convenient to number the vertices of the $n$-gon in the definition of
$\Delta^\star_{n,k}$ from $0$ to $n-1$ and the vertices of the $(n-1)$-gon 
in the definition of $\Delta_{n-1,k}$ by $1$ to $n-1$.

By a subtle observation of Jonsson \cite{Jopre} there exist sets of
diagonals $B,B_1$ in the $n$-gon resp. $(n-1)$-gon such that
\begin{eqnarray}
\link_{\Delta^\ast_{n,k}}B\cap\Omega_{n,k}&\cong&\link_{\Delta^\ast_{n-1,k}}B_1\cap\Omega_{n-1,k}.\label{linkiso}
\end{eqnarray}
Recall that the link of a complex $\Delta$ with vertex-set $V$ with respect
to a face $\sigma\in\Delta$ is defined as
$$\link_{\Delta}(\sigma):=\{\tau\in \Delta\big|\tau\cup\sigma\in
\Delta,\tau\cap\sigma=\emptyset\}.$$

The definition of the sets $B,B_1$ is easier with the help of the set
$\gamma_{n,k}:=\Gamma_{n,k}\setminus\Gamma_{n,k+1}$.

\begin{lemma}[\cite{Jopre,Jo}]\label{Jaklemma1} For every facet $F$ in
$\Delta^\ast_{n,k}$ exists a unique set $B:=B_1\cup B_0$ with
$$B_1:=\{1b_1,2b_2,\dots,kb_{k}\}$$
and $$B_0:=\{0b_1,1b_2,\dots,k-1b_{k}\},$$ such that 
$$k+1\leq b_1<b_2<\dots<b_k\leq (n-1) \textrm{ and } B\subset F\cup\gamma_{n,k}.$$

If we additionally define
\begin{eqnarray*}Z_B&:=&\{ij\big|b_i<j<b_{i+1}, i=0,\dots, k\}\\
K_B&:=&\{ij\big| k+1\leq i<b_1,b_k<j\leq n-1\}\end{eqnarray*}
we get that
\begin{eqnarray*}
Z_B\cap F&=&\emptyset\\
K_B\cap F&=&\emptyset.\end{eqnarray*}
\end{lemma}

The isomorphism in (\ref{linkiso}) is defined with the help of a partition of the
vertex-sets as follows:

\begin{lemma}[\cite{Jopre,Jo}]\label{Jaklemma2}
For a given facet $F$ of $\Delta^\ast_{n,k}$ let $B_0,B_1,K_B,Z_B$ be defined as
in Lemma \ref{Jaklemma1}. We get the partitions
\begin{eqnarray*}
\Omega_{n,k}&=&B\cup Z_B\cup K_B\cup S_0\cup S_1\cup S_2\\
\Omega_{n-1,k}&=&B_1\cup K_B\cup S'_0\cup S_1\cup S_2\\
V_{n,k}&=&S_0\cup S_1\cup S_2\\
V_{n-1,k}&=&S'_0\cup S_1\cup S_2.
\end{eqnarray*}
where the latter sets are the vertex-sets
$V_{n,k}$ of $\link_{\Delta^\ast_{n,k}}(B\cap\Omega_{n,k})$ and
$V_{n-1,k}$ of  $\link_{\Delta^\ast_{n-1,k}}(B_1\cap\Omega_{n-1,k})$
by defining
\begin{eqnarray*}S_0&:=&\{(i-1)j\big|i\in [1,k],b_i<j\leq n-k-2+i\}\\
S'_0&:=&\{ij\big|i\in [1,k],b_i<j\leq n-k-2+i\}\\
S_1&:=&\{ij\big|1\leq i\leq k,i+k+1\leq j <b_i\}\\
S_2&:=&\{ij\big|k+1\leq i,j\leq n-1\}\cap \Omega_{n,k})\backslash
K_B.\end{eqnarray*}
\end{lemma}

\begin{lemma}[\cite{Jopre,Jo}]\label{Jabb}
The mapping $$\varphi_B:S_0\cup S_1\cup S_2\cup B \longrightarrow S'_0\cup
S_1\cup S_2\cup B_1$$
$$ij\mapsto\begin{cases}(i+1)j&\textrm{ if }ij\in B_0\cup S_0\\ij&\textrm{
if }ij\in B_1\cup S_1\cup S_2\end{cases}$$
induces an isomorphism
$$\link_{\Delta^\ast_{n,k}}B\cap\Omega_{n,k}\cong\link_{\Delta^\ast_{n-1,k}}B_1\cap\Omega_{n-1,k}.$$
\end{lemma}

An immediate consequence of the foregoing lemmas is
\begin{lemma}\label{Mengenlemma}
For a given facet $F$ and its unique set $B$,
the vertex-set $V^{symm}_{2n,k}$ of
$\link_{\Delta^{symm}_{2n,k}}(B\cup\bar
B)\cap\Omega_{2n,k}))$ can be written as
\begin{eqnarray*}
V^{symm}_{2n,k}&=&(S_0\cup S_1\cup S_2)\cap(\bar S_0\cup \bar S_1 \cup \bar
S_2).
\end{eqnarray*}
\end{lemma}

Again it turns out to be convenient to let $\Delta^{\ast\text{symm}}_{2n,k}$
be coming from a $2n$-gon numbered from $0$ to $2n-1$ and
 $\Delta^{\ast\text{symm}}_{2n-2,k}$
be defined on the $2n-2$-gon with numbering of the vertices from $1$ to
$n-1$ and from $n+1$ to $2n-1$.

\begin{figure}
\includegraphics[scale=0.45]{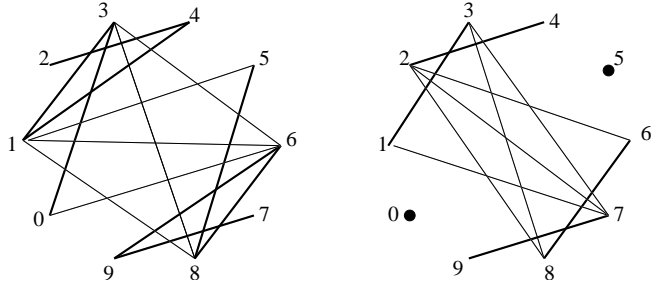}\
\caption{Action of $\varphi$ on a facet of $\Delta^{\ast\text{symm}}_{5,2}$,
elements of $B\cup\bar B$ printed bold.}\label{phiinaction}
\end{figure}

We define the effect of $$\varphi:V^{symm}_{2n,k}\cup B\cup \bar
B\rightarrow \Omega_{2n-2,k}$$ on the diagonal $ij$ with $i<j$ by the
following table, each cell corresponding to the intersection of the sets in
the respective row and column, this intersection being empty for cells
marked with $\Phi$.
The mapping $\varphi$ can be seen in action in Figure \ref{phiinaction}.

\begin{center}
\begin{figure}
\begin{tabular}{c|ccc|cc}\label{table1}
 &$S_0$&$S_1$&$S_2$&$B_0$&$B_1$\\
 \hline
 &&&$\,$ \\
 $\bar S_0$&$(i+1)(j+1)$ &$i(j+1)$&$i(j+1)$&$(i+1)(j+1)$&$i(j+1)$\\
 $\bar S_1$&$(i+1)j$&$ij$&$ij$&$(i+1)j$&$ij$\\
 $\bar S_2$&$(i+1)j$&$ij$&$ij$&$(i+1)j$&$ij$\\
 &&&$\,$ \\
 \hline
 &&&$\,$ \\
 $\bar B_0$&$(i+1)(j+1)$&$i(j+1)$&$i(j+1)$&$(i+1)(j+1)$&$\Phi$\\
 $\bar B_1$&$(i+1)j$&$ij$&$ij$&$\Phi$&$ij$\\

\end{tabular}
\caption{Definition of $\varphi$.}

\end{figure}
\end{center}

\begin{lemma}\indent\label{phiproperty}

\begin{enumerate}
\item \label{i} The mapping $\varphi$ preserves rotational symmetry, i.e. $\varphi(\bar
e)=\overline{\varphi(e)}$ for all diagonals $e$ in the domain of
$\varphi$.

\item \label {ii} Each diameter in $\Omega_{2n-2,k}$ is the image of a
diameter in $V_{2n,k}^{symm}\cup B\cup\bar B$.

\item \label{iii} For a fixed facet $F$ of $\Delta_{2n,k}^{\ast\text{symm}}$, each
diameter in $\varphi(F)$ is the image of a diameter in $F$.

\item \label{iiii} If a subset $\sigma\subset V_{2n-2,k}^{symm}\cup B\cup\bar
B$ does not contain a $(k+1)$-crossing, this is also true for $\varphi(\sigma)$.

\item \label{iiiii} If a subset $\sigma\subset \Omega_{2n-2,k}$
does not contain a $(k+1)$-crossing, this is also true for $\varphi^{-1}(\sigma)$.

\end{enumerate}
\end{lemma}

\begin{proof}\indent
\begin{itemize}
\item[(\ref{i})] This is clear from the symmetry of Table \ref{table1}.

\item[(\ref{ii})] Let $d=i(i+n)$ be a diameter in $\Omega_{2n-2,k}$ and
$d'=(i-1)(n+i-1)\in\Omega_{2n,k}$.
We need to discriminate several cases:
\begin{itemize}

    \item[(A)] $i\notin\{1,\dots,k\}$. Then we have $d\in S_2\cap\bar S_2$ and
                           $\varphi(d)=d$.

     \item[(B)] $i\in\{1,\dots,k\}$

       \begin{itemize}

        \item[(B1)] If $b_i<n+i-1$ then $d'\in
        S_0\cap\bar S_0$ and $\varphi(d')=d$.

    \item[(B2)] If $b_i=n+i-1$ then $d'\in B_0\cap \bar B_0$,
                          and $\varphi(d')=d$.

\item[(B3)] If $b_i=n+i$, then $d\in B_1\cap \bar B_1$ and
$\varphi(d)=d$.

\item[(B4)] If $b_i>n+i$, then $d\in S_1\cap \bar S_1$ and
$\varphi(d)=d$.

\end{itemize}

\end{itemize}

\item[(\ref{iii})] Let $F$ be a facet in $\Delta_{2n,k}^{symm}$ and $d=i(i+n)$ a diameter in $\varphi(F)$.
Whenever we have $\varphi(e)=d$ for a non-diameter $e\in\Omega_{2n,k}$, then
we also have $\varphi(\bar e)=d$ and we know $e\in S_0\cup B_0\setminus (\bar S_0\cup \bar B_0)$.
We know there is an $i$ such that $e,\bar
e$ are of the form $e=(i-1)(i+n),\bar e=i(i+n-1)$.
Any diagonal of the $2n$-gon, except for
$d':=(i-1)(n+i-1)$, which crosses $d$, crosses at least
one of the diagonals $e,\bar e$. Since $d,d'$ do not cross $e,\bar e$, one
of $d$ and $d'$ is contained in $F$.

If $e\in B_0$, we have $i-1\in[0,k-1]$ and $n+i=b_{i}$. Then $d\in B_1\cap
\bar B_1$ and as a consequence $d\in F$ and $\varphi(d)=d$.

If $e\in S_0$, we have $i-1\in[0,k-1]$ and $b_{i}<n+i$. But $b_i$ can not be
smaller than $n+i-1$: Since $\bar e\notin S_0\cup B_0$, we have $n+i-1<b_{i+1}$, but this means
$\bar e\in Z_B\cap F=\emptyset$ (see Lemma \ref{Jaklemma1}).

As a consequence we have $b_i=n+i-1$, meaning that $d'\in B_0\cap \bar B_0$
and by this $d'\in F$ and $\varphi(d')=d$.

\item[(\ref{iiii})]
Let $\sigma\subset V_{2n-2,k}^{symm}\cup B\cup\bar B$ be
a maximal subset not containing a $(k+1)$-crossing. Assume
$\tau=\varphi(\sigma)$ comprises a $(k+1)$-crossing $E$, which we choose in
a fashion such that the number of elements in $\varphi(B_0)$ is maximal compared
to other $(k+1)$-crossings in $\tau$.

Since $\varphi(B_0)=\varphi(B_1)$ as one easily shows, we get
$\varphi(\sigma\setminus B_0)=\varphi(\sigma)$. As a consequence we have two
non-crossing segments $ix,jy\in \sigma\setminus
B_0$ such that $\varphi(ix),\varphi(jy)\in E$. W.l.o.g. we can conclude that
$i=j\leq k-1, ix\in S_0\setminus B_0$ and $iy\notin S_0\cup B_0$. This
implies $y\leq b_{i}$ since otherwise $iy\in Z_B$ and $x>b_{i+1}$ as well as
$\varphi(ix)\in\{(i+1)x,(i+1)(x+1)\}$.

We define sets $E',E''$ by
\begin{eqnarray*}E'&:=&\left(E\setminus \varphi(ix)\right)\cup (i+1)b_{i+1}\\
E''&:=&\left(E\setminus\varphi(ix)\right)\cup (i+1)(b_{i+1}+1).\end{eqnarray*}

If $x>b_{i+1}+1$, both $E$ and $E'$ are $(k+1)$-crossings. But at least one
of $E'$ and $E''$ contains an element of $\varphi(B)$ which is not in $E$ in
contradiction to the choice of $E$.

That leaves us with $b_{i+1}+1$ as the only possible value for $x$.

In this case $E'$ is a $(k+1)$-crossing different from $E$, consequently
$(i+1)b_{i+1}$ must not be an element of $\varphi(B)$.

The set $E''$ is a new $(k+1)$-crossing with more elements from $\varphi(B)$
if $\varphi(ix)=(i+1)(x+1)$, because else we have $E''=E$.

Thus we need to have
\begin{enumerate}

\item $(i+1)b_{i+1}\notin \varphi(B)$\label{Proba}
\item $\varphi(ix)=(i+1)x$\label{Probb}
\item $x=b_{i+1}+1$.\label{Probc}

\end{enumerate}

From (\ref{Proba}) we can conclude that $ib_{i+1}\in B_0\cap (\bar B_0\cup\bar
S_0)$. A consequence of (\ref{Probb}) and (\ref{Probc}) is $ix\in S_0\setminus(\bar
B_0\cup\bar S_0)$.

All in all this leads to $n\leq b_{i+1}\leq n+k-1, i\geq b_{b_{i+1}-n}$ and
$b_{i+1}+1>n+k-1$. We find that $b_{i+1}=n+k-1$ and thus $i\geq b_{k-1}\geq
k+1$ what finally yields a contradiction.

\item[(\ref{iiiii})] First we show that two crossing diagonals still cross after the
application of $\varphi$:
Let $ix$ and $jy$ be crossing diagonals, while $\varphi(ix),\varphi(jy)$
do not cross. W.l.o.g. we can assume
$j=i+1,\,ix\in S_0\cup B_0,(i+1)y\in B_1\cup S_1\cup S_2, y>x.$
As a consequence we get  $x\geq b_{i+1}$ and, for the case $(i+1)y\in B_1\cup S_1$,
that $y\leq b_{i+1}$, contradicting $y>x$.
If $(i+1)y$ was an element of $S_2$, we had $i+1=k$ and from $y>x\geq b_k$
we conclude that $(i+1)y\in K_B\not\subset S_2$.

Now we assume that in contrast to the assertion a subset
$\sigma\subset \varphi(\Omega_{2n,k})$ does not contain a
$(k+1)$-crossing and $\varphi^{-1}(\sigma)$ does. Then with the result above
$\varphi(\varphi^{-1}(\sigma))=\sigma$ contains a $(k+1)$-crossing as well,
contradicting the assumption.

\end{itemize}\end{proof}

With Lemma \ref{phiproperty} we are ready for the proof of Theorem
\ref{diameters}:

\begin{proof}[Proof of Theorem \ref{diameters}] For arbitrary $k$ we proceed by induction over $n$, the basecase being
$n=2k+2$. Here $\Omega_{2k+2,k}$ exclusively consists of diameters, all
mutually intersecting, so that each
facet contains exactly $k$ of them.
Now for $n>2k+2$, choose a facet $F$ from $\Delta^{\ast\text{symm}}_{2n,k}$. Let $B$ as in
Lemma \ref{Mengenlemma} and $\varphi$ be defined on the corresponding
partition. Then we know by Lemma \ref{phiproperty} (\ref{iiii}), that
$\varphi(F)$ does not contain a $(k+1)$-crossing.
Adding another diameter $d$ to $\varphi(F)$ would produce a $(k+1)$-crossing,
since otherwise according to Lemma \ref{phiproperty}, (\ref{ii}), $\varphi^{-1}(d)$
exists and $\varphi^{-1}(\varphi(F)\cup\{d\})$ would be a proper superset of
$F$ and be free of $(k+1)$-crossings, contradicting the facet-property of
$F$.
According to our assumption, $\varphi(F)$ contains exactly $k$ diameters,
each of them is the image of a diameter in $F$ as stated in ($\ref{iii}$) of
the foregoing lemma. This means that $F$ comprises at least $k$ diameters
and since all diameters are mutually intersecting there are exactly $k$ diameters.
\end{proof}

\section{Proof of Theorem \ref{homsphere}} \label{proofhomsphere}

For a simplicial complex $\Delta$ and a finite group $G$ we call $\Delta$ a
$G$-complex if $G$ acts simplicially on the vertex-set of
$\Delta$, i.e. $g\sigma:=\{ gv \big|v\in\sigma\}\in \Delta$ for all
$\sigma\in\Delta$ and $g \in G$. We call a simplicial $G$-complex regular if
for each subgroup $U\leq G$ of $G$ and any choice of elements $g_0,\dots,g_n\in U$ 
we have that if $\{v_0,\dots,v_n\}$ and $\{g_0v_0,\dots,g_nv_n\}$ are both 
simplices in $\Delta$, there exists an element $g$ in $U$ such that $gv_i=g_iv_i$
for all $0\leq i\leq n$. 

Recall that the barycentric subdivision of a simplicial complex $\Delta$
is a simplicial complex on vertex-set $\Delta \setminus \{ \emptyset \}$ 
whose simplices are the subsets of $\Delta \setminus \{ \emptyset \}$
that are totally ordered with respect to inclusion.
Clearly, if $\Delta$ is a $G$-complex then $\sd(\Delta)$
is a $G$-complex as well.  We will make use of the following theorem.

\begin{theorem}[\cite{Br}]\label{vorbred} Let $\Delta$ be a simplicial
$G$-complex. If for all $g\in G,\,\sigma\in \Delta$ and all $\sigma \in 
\Delta$ we have $gv=v$ for all vertices $v\in\sigma\cap g(\sigma)$, then 
the barycentric subdivision $\sd(\Delta)$ is a regular $G$-complex.
\end{theorem}

For a $G$-complex $\Delta$ its fix-complex $\Delta^G$ is the simplicial complex consisting of
those simplices $\sigma\in\Delta$ that are elementwise fixed by $G$.
For a regular $G$-complex its fix-complex sometimes inherits topological properties from 
the complex.

\begin{theorem}[\cite{Br}] \label{bredon51} Let $G = \ZZ_p$ be the cyclic group of prime order 
$p$.  If $\Delta$ is a $d$-dimensional regular simplicial $G$-complex such that
$H_i(\Delta,\ZZ_p)\cong H_i(S^d,\ZZ_p)$ for all $i\leq d$ then there
is an $\ell\leq d$ such that $H_i(\Delta^G,\ZZ_p)\cong
H_i(S^\ell,\ZZ_p)$ for all $i\leq d$.
\end{theorem}

Note, that in contrast to our definition, in \cite{Br} a simplicial complex with the homological 
properties required in Theorem \ref{bredon51} is called a $\mot p$-homology-sphere.

We will apply Theorem \cite{Br} to a suitable subdivison of $\Delta^\ast_{2n,k}$. 
Since by Theorem \ref{sphere} a geometric realization of $\Delta^\ast_{2n,k}$ 
is a sphere, it follows immediately that 
is a $\mot p$-homology-sphere for all primes $p$.

\begin{figure}
\includegraphics[scale=0.3]{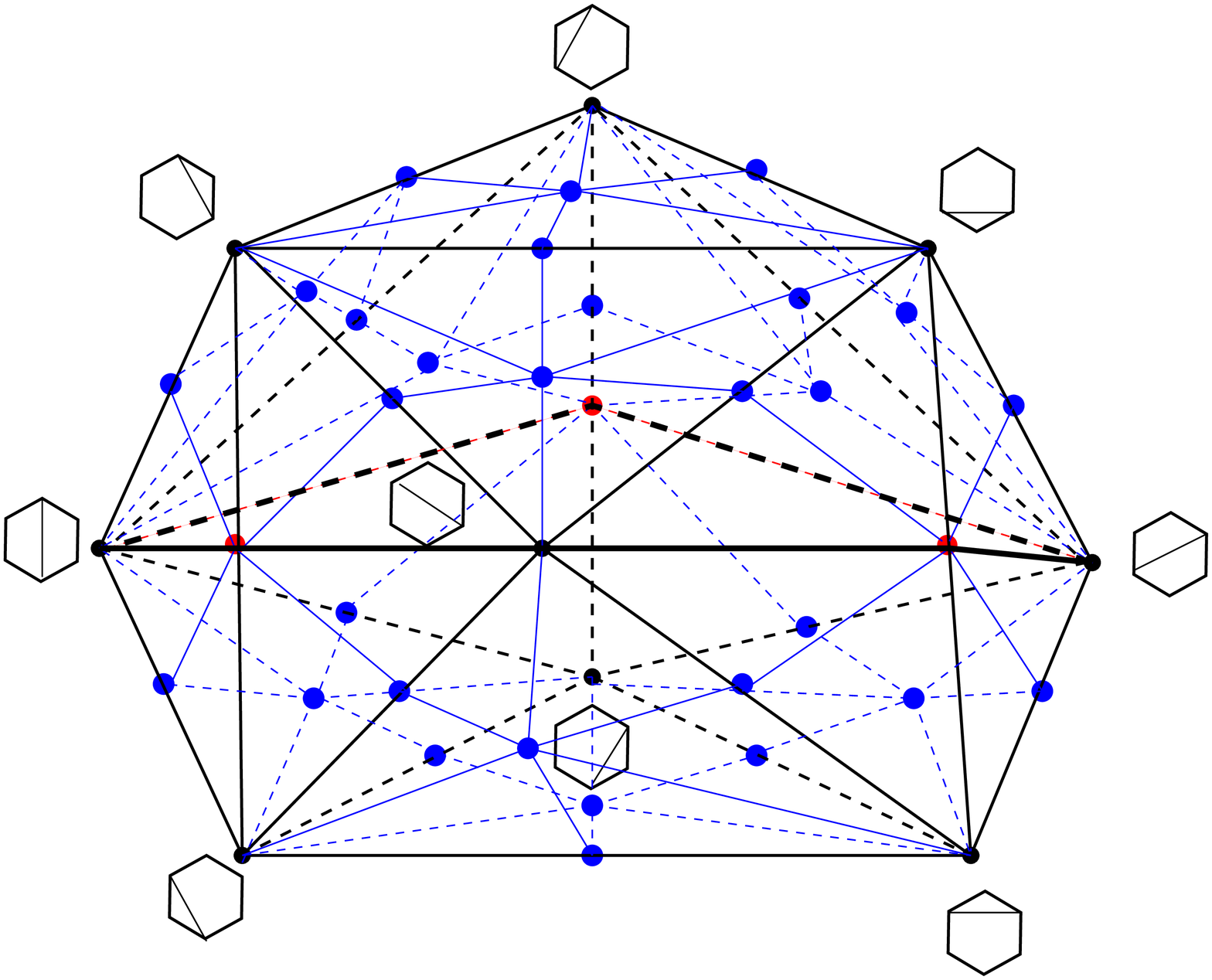}\
\caption{Construction in the proof of Theorem \ref{homsphere}}\label{assoziaedermit}
\end{figure}

For the proof of Theorem \ref{homsphere} we construct a subdivision 
$S_{2n,k}$ of the complex $\Delta^\ast_{2n,k}$, such that the group 
$\ZZ_2$
act simplicially and regularly and $\D_{n,k}^\ast$ is isomorphic to the fix-complex $S_{2n,k}^{\ZZ_2}$.
The assertion then is a consequence of Theorem \ref{bredon51} and Corollary \ref{dimtypb}.
Figure \ref{assoziaedermit} shows $\Delta^\ast_{6,1}$ and the construction of
$T_{6,1}$. The bold facets are the faces of the fix-complex
and correspond to faces of $\D^\ast_{3,1}$.

The complex $\Delta^\ast_{2n,k}$ carries a natural $\ZZ_2$-action which is induced
by sending a diagonal $d$ to its image $\bar{d}$ under $180^\circ$-rotation. From now on 
we identify $\Delta^\ast_{2n,k}$ with its geometric realization and also consider the 
vertices as points in some $\RR^m$. 
We set $$\CR_{2n,k} := \{\frac{1}{2}(d+\bar d)\big|d\in\Omega_{2n,k} ,d\not=\bar d\}, 
\Omega_{2n,k}':=\Omega_{2n,k} \cup \CR_{2n,k}.$$
The set $\Omega_{2n,k}'$ will serve as the ground set of our subdivision.
For $\sigma \in \Delta^\ast_{2n,k}$ we let $\sigma^{symm}:=\{d\in\Omega_{2n,k} \big|d\in\sigma\textrm{ and }\bar d\in
\sigma,d\not=\bar d\}$
be the symmetric part of $\sigma$ without the diameters
and set $2l:=\sharp\sigma^{symm}$. We choose an $l$-element subset $D_\sigma:=\{d_1,\dots,d_l\}$ from
$\sigma^{symm}$ such that $d_i\not=\bar d_j$ for all $1\leq i<j\leq l$.
For each $d_i$ and each $\varepsilon\in \{0,1\}^l$ we set 
$$d_i^{\varepsilon_i}:=\begin{cases}d_i\textrm{ if
}\varepsilon_i=1\\ \bar d_i\textrm{ if }\varepsilon_i=0\end{cases}$$
and $D_\sigma^\varepsilon:=\{d_i^{\varepsilon_i},i=1,\dots,l\}$. Finally we
define for $\varepsilon\in \{0,1\}^l$ the simplex
$$\sigma_\varepsilon:=\left( \sigma\setminus\sigma^{symm} \right) \cup\{\frac{1}{2}(d+\bar d)\big|d\in D_\sigma\}\cup D_\sigma^\varepsilon.$$
We let $T_{2n,k}$ be the simplicial complex which is generated by all $\sigma_\varepsilon$ for $\sigma \in \Delta^\ast_{2n,k}$
and $\varepsilon\in \{0,1\}^l$.  

By 
$\conv(\{D_\sigma^\varepsilon\big|\varepsilon\in\{0,1\}^l\})=\conv(\sigma^{symm})$,
it follows that $T_{2n,k}$ is a subdivision of $\Delta^\ast_{2n,k}$. 
Since as mentioned before $\Delta^\ast_{2n,k}$ is a 
a $\mot 2$-homology-sphere the same is true for $T_{2n,k}$..

Clearly, $\ZZ_2$ acts simplicially on $T_{2n,k}$.

Since each $\sigma_\epsilon$ contains at most one of $d$ and $\bar d$ for any 
$d \in \Omega_{2n,k}$, we get 
for all $\sigma \in T_{2n,k}$ and all $g \in \ZZ_2$ that 
\begin{eqnarray}\label{fix}gv=v\text{ for all }v\in
g\sigma\cap\sigma.
\end{eqnarray}

Now by a suitable subdivision we give $T_{2n,k}$ the structure of
a regular $\ZZ_2$-complex. 

We replace all $\sigma\in T_{2n,k}$ by the join $\sigma_1 * \sigma_2$, where
$\sigma_1$ is the barycentric subdivision of $\sigma \cap \Omega_{2n,k}$ and
$\sigma_2 = \sigma \cap \CR_{2n,k}$. Note, that here we adopt the convention
that the join of a simplex with the empty set is the simplex.
We write $S_{2n,k}$ for the resulting simplicial complex. By construction,
$S_{2n,k}$ is a subdivision of $T_{2n,k}$ and therefore of $\Delta_{2n,k}^\ast$.

To show regularity it suffices to show the following. If 
$\sigma:=\{v_0,\dots,v_n\}$ and $\tau:=\{v_0, \ldots, v_r, \bar{v_{r+1}}, \ldots, \bar{v_n} \}$
be two simplices of $S_{2n,k}$ then either $\bar{v_{i}} = v_i$ for $r+1 \leq i \leq n$ (resp. 
$\sigma = \tau$) or $\bar{v_{i}} = v_i$ for $0 \leq i \leq r$ (resp. $\sigma = \bar{\tau}$). 

Since $\CR_{2n,k}$ is pointwise fixed by $\ZZ_2$ we can ignore the vertices from 
$\CR_{2n,k}$ and assume that $\sigma, \tau$ are simplices in the barycentric
subdivision of $T_{2n,k}' := \{ \rho \in T_{2n,k}~|~\rho \subseteq \Omega_{2n,k}\}$.
By (\ref{fix}) the complex $T_{2n,k}$ and therefore $T_{2n,k}'$ satisfy the assumption of Theorem \ref{vorbred}.
The theorem implies that the barycentric subdivision of $T_{2n,k}'$ is a regular
$\ZZ_2$-complex. From this it follows that either $\sigma = \tau$ or $\sigma = \bar{\tau}$. 

By construction, 
$$\phi:\{d,\bar d\} \mapsto \frac{1}{2}(d+\bar d)$$
induces an isomorphism between $\D^\ast_{n,k}$ and
the fixcomplex $S_{2n,k}^{\ZZ_2}$. With Theorem \ref{bredon51} we get
that $$H_i(\mathcal D^\ast_{n,k},\ZZ_2)\cong H_i(S^{k(n-k)-1},\ZZ_2).$$

Now let $\tau$ be a face of $\D^\ast_{n,k}$. Since $\D^\ast_{n,k}$ is pure of dimension
$k(n-k)-1$, we have $\dim\link_{\D^\ast_{n,k}}(\tau)=k(n-k)-1-\dim\tau-1$.

The isomorphism from above yields
$$\link_{\D^\ast_{n,k}}(\tau)\cong \link_{S_{2n,k}^{\ZZ_2}}(\phi(\tau))$$
and since  $S_{2n,k}^{\ZZ_2}$ is a $\mot 2$-homology-sphere, the
assertion follows.

\section{Proof of Theorem \ref{Stan} and Proposition \ref{polytopaltypb}} \label{proofStan}.

  The statements of Theorem \ref{Stan} and Proposition \ref{polytopaltypb}
  seem to be unrelated. Nevertheless, we provide their proofs in a joint
  section, since the proof of Proposition \ref{polytopaltypb} makes use of a
  bijection established in the proof of Theorem \ref{Stan}.

  \begin{proof}[Proof of Theorem \ref{Stan}] 
    Let $V=[n]\times [n]$. It is easy though tedious to check that the mappings
    $\Psi:\mathcal F_n\rightarrow V$ and $\Phi:V\rightarrow
    \mathcal F_n$ given by 
    \begin{eqnarray*}
       \Psi:\quad\{\{a_1<b_1\},\{a_2<b_2\}\}\mapsto \begin{cases}((a_i+1)\mot n,b_i\mot
       n)&\textrm{ if }b_i-a_i\leq n\\((b_i+1)\mot n,a_i\mot n)&\textrm{ if
       }b_i-a_i>n,\end{cases}\\
       \Phi:\quad(a,b)\mapsto \begin{cases}\{\{(a-1)\mot 2n,b\},\{a-1+n,b+n\}\}&\textrm{
       if }a\leq b \\
       \{\{a-1+n,b\},\{a-1,b+n\}\}&\textrm{ if }a>b\end{cases}
     \end{eqnarray*}

    are well-defined mutually inverse bijections.

    Furthermore, $\Psi$ maps the set $F_{n,k}$ to $[n]\times [n]\backslash
    \{(i,j)~\big|~j=(i+l)\mot n\textrm{ for one }l<k\}$ and for a $(k+1)$ crossing
    $K\subset \mathcal F_{n,k}$ there are $(k+1)$-sets $A,B\subset [n]$ such that 
    $\Psi(K)=N(A,B)$. Conversely, any subset of $V$ of the type $N(A,B)$ is
    mapped on a $(k+1)$-crossing by $\Phi$. This already proves Theorem \ref{Stan}.
  \end{proof}

  \begin{proof}[Proof of Proposition \ref{polytopaltypb}] 
     If $k =1$ then $D_{n,k}^\ast$ is the cyclohedron, which is well known 
     (see for example \cite{Si}) to be the boundary complex of a polytope.
     In the case $k = n-1$ the complex $D_{n,k}^\ast$ is the boundary complex of
     the simplex on the set of diameters. it remains
     to prove the proposition for $n\geq 3$ the complex $\D^\ast_{n,n-2}$ is isomorphic 
     to the boundary-complex of the $(2n-4)$-dimensional cyclic polytope with $2n$ 
     vertices $\mathcal C_{2n-4}(2n)$.
     We use the following characterization of the boundary-complex by Gale:
     Let $C_{2n}$ the cycle-graph with $2n$-vertices numbered in clockwise order
     and for $\sigma\subset [2n]$ let $C_{2n}(\sigma)$ be the subgraph induced by
     the vertices in $\sigma$.

     \begin{theorem}[see \cite{Gr}] \label{Galeeven} 
       Identifying each vertex
       $x_i=(i,i^2,i^3,\dots,i^d)\in\RR^d$ of $\mathcal C_{2n-4}(2n)$ with
       vertex $i$ of the graph $C_{2n}$ for
       $i=1,\dots,2n$ yields:
       The subset $\sigma\subset [n]$
       is a face of the boundary-complex of $\mathcal C_d(2n)$ if and only if 
       $|\sigma|+\omega(\sigma)\leq d,$ where $\omega(\sigma)$ is the number of
       odd-sized connected components of $C_{2n}(\sigma)$.
     \end{theorem}

     The correspondence from Theorem \ref{Galeeven} translates into the following classification
     of the minimal nonfaces of $\mathcal C_{2n-4}(2n)$. 

     \begin{corollary}\label{Nichtseitenzyklisch}
       The minimal nonfaces of $\mathcal C_{2n-4}(2n)$ are in bijection with 
       the $(n-1)$-subsets $\eta\subset[2n]$, for which $C_{2n}(\eta)$ consists of
       exactly $(n-1)$ connected components.
     \end{corollary}

     As noted in the paragraph above, we view $\D^\ast_{n,n-2}$ via the
     bijection $\Psi$ from the proof of Theorem \ref{Stan} as a simplicial
     complex on vertex-set $[n]\times [n]\backslash \{(i,j)~\big|~j=(i+\ell)\mot
     n\textrm{ for
     a }\ell<n-2\}$ with nonfaces $N(A,B)$ for the $(n-1)$-element subsets $A$, $B$ of
     $[n]$, .
     If seen as a set of entries in a $n\times n$-matrix this vertex set 
     is a set of matrix elements such that there are $2$ entries in each row and column.

     We number the vertices such that for each
     $i\in\{1,\dots,2n\}$ the vertex with number $i$ shares a row or a column
     with the vertices numbered by $i-1,i+1 (\mot 2n)$. Recall that $k\mot
     k=k$ (See Figure \ref{nurzwei} for an example.)

     \begin{figure}
       \begin{center}
         $\left(\begin{array}{cccccc}
         \star &\star & \star &\star&2 &1 \\
         11 &\star &\star &\star&\star &12 \\
         10 & 9 &\star &\star&\star&\star \\
         \star &8 & 7 &\star&\star &\star \\
         \star &\star & 6 &5&\star &\star \\
         \star &\star & \star & 4 &3 &\star \end{array}\right)$
       \end{center}
       \caption{Numbering of the vertices of $\mathcal D^\ast_{6,4}$ as in the proof of Proposition 
          \ref{polytopaltypb}. }\label{nurzwei}
     \end{figure}

     We now show that each numbering $\varphi$ of the vertices with this property identifies 
     nonfaces of $\D^\ast_{n,n-2}$ with nonfaces of $\mathcal C_{2n-4}(2n)$ and vice versa.

     For all $A,B\subset [n]$ with $|A|=|B|=n-1$ the set
     $N(A,B)$ contains exactly $n-1$ elements and in each row and column of $[n]\times [n]$
     there is at most one element of $N(A,B)$. Hence the subgraph $C_{2n}(\varphi(N(A,B))$
     does not contain neighbouring vertices of $C_{2n}$ and thus induces exactly $(n-1)$ connected components.

     Conversely, if we have a $(n-1)$-subset $\tau$ of the vertex-set of $C_{2n}$
     such that $C_{2n}(\tau)$ has exactly $(n-1)$ connected components, then
     elements from $\varphi^{-1}(\tau)$ are spread over $(n-1)$ different rows
     and columns. Let $A=\{a_1,\dots,a_{n-1}\}$ be the set of row-indices and
     $B=\{b_1,\dots,b_{n-1}\}$ the set of column-indices.

     Both sets $\varphi^{-1}(\tau)$ and $N(A,B)$ have the property that no two
     of their elements share a row or a column. Together with the fact
     that in each row and column of the matrix $[n]\times [n]$ there are only two
     vertices of $\D_n^{n-2}$ we get that the two sets coincide whenever they have
     nonempty intersection.

     Define $i_A,i_B \leq n$ to be the unique elements of $[n] \setminus A$ and $[n] \setminus B$.
     First, we examine the set $\varphi^{-1}(\tau)$. If $i_A\not=(i_B+2)\mot n$
     we have $((i_B+2)\mot n,(i_B+1)\mot n)\in\varphi^{-1}(\tau)$, since it is
     the unique entry of the matrix in row  $(i_B+2)\mot n$.
     If $i_A=(i_B+2)\mot n$, we have that $((i_B+1)\mot
     n,(i_B-1)\mot n)$ is the unique entry in row $(i_B+1)\mot n$ and thus
     an element of $\varphi^{-1}(\tau)$.

     On the other hand the set $N(A,B)$ is of the form $\{(a_{(j+\ell)\mot
     n-1},b_j)~|~j=1,\dots,n-1\}$ for some $0\leq \ell \leq 2$.

     \begin{itemize}
       \item[$\ell = 0$:]
         We have $\ell=0$ if and only if $i_A=1$ and $i_B=n$, since only in this case
         $a_j=j+1>j=b_j$ holds for all $j=1,\dots,n-1$.
         As a consequence we get that  $((i_B+2)\mot n,(i_B+1)\mot n)=(2,1)=(a_1,b_1)\in
         N(A,B)\cap\varphi^{-1}(\tau)$.
      \item[$\ell =2$:] We have $\ell = 2$ if and only if $i_B\leq i_A-2$, since only in 
         this case there is a $j\in\{i_B,\dots,i_A-2\}\subset\{1,\dots,n-2\}$
         such that $a_{j+1}\leq b_j$ and we always have
         $a_{i+2}\geq i+2>i+1\geq b_i$ for all $i\in\{1,\dots,n-3\}$.

         For $i_B<i_A-2$ we have $a_{i_B+2}=i_B+2$ und $b_{i_B}=i_B+1$ and as a
         consequence $(i_B+2,i_B+1)\in N(A,B)\cap \varphi^{-1}(\tau)$.
         For $i_B=i_A-2$ we can conclude that $a_{i_B+1}=i_B+1$ and
         $b_{i_B-1}=i_B-1$. Thus $(i_B+1,i_B-1)\in N(A,B)\cap
         \varphi^{-1}(\tau)$.

       \item[$\ell = 1$:] By the previous argumentation we have $\ell=1$ for the remaining cases, that
         is $i_B>i_A-2$ and $(i_A,i_B)\not=(1,n)$.
         If $i_B<n-1$ (and consequently $i_A\not=(i_B+2)\mot n$) we have
         $(i_B+2,i_B+1)=(a_{i_B+1},b_{i_B})\in N(A,B)\cap \varphi^{-1}(\tau).$
         If $i_B=n-1$ and $i_A\not=(i_B+2)\mot n=1$ we have
         $(1,n)=(a_{1},b_{n-1})\in N(A,B)\cap \varphi^{-1}(\tau).$

         For $i_B=n$, and $i_A\not=(i_B+2)\mot n=2$ we get
         $$(2,1)=(a_{2},b_{1})\in N(A,B)\cap \varphi^{-1}(\tau).$$
         Whenever $i_B=n-1$, and $i_A=(i_B+2)\mot n=1$ we conclude
         $$(n,n-2)=(a_{n-1},b_{n-2})\in N(A,B)\cap \varphi^{-1}(\tau).$$
         Finally, if  $i_B=n$, and $i_A=(i_B+2)\mot n=2$ we know
         $$(1,n-1)=(a_{1},b_{n-1})\in N(A,B)\cap \varphi^{-1}(\tau).$$
     \end{itemize}
   \end{proof}

\section{Construction of the term-order}\label{definition}

  \begin{defn}\label{varordnung}
    Let
    \begin{eqnarray*} \varphi:[n]\times [n]&\rightarrow& \NN\\
      (i,j)&\mapsto&\left[(2-i)\cdot n+(j-1)\cdot (n-1)-1\mot n^2\right]+1,
    \end{eqnarray*}
    and order the entries in $X=(x_{ij})_{1\leq i,j\leq n}$ according to
    $$x_{ij}\preceq x_{kl}:\Leftrightarrow \varphi(i,j) \leq \varphi(k,l).$$

    \begin{figure}
      \begin{center}
        $\left(\begin{array}{ccccc}
        5 &    9 &   13 &   17 &   21 \\
        25  &   4 &   8 &   12 &   16 \\
        20  &  24 &    3 &    7 &   11 \\
        15  &  19 &   23 &    2 &    6  \\
        10  &  14 &   18 &   22 &    1\end{array}\right)$
      \end{center}
      \caption{The mapping $\varphi$ for $n=5$}
    \end{figure}
  \end{defn}

  It is easily checked that $\phi$ is a bijection from $[n] \times [n]$ to $[n^2]$.

  We now define an order on the set of monomials, which will be the
  graded, reverse-lexicographic continuation of $\preceq$. For that we denote each
  monomial $t=\prod_{(i,j)\in[n]^2}x_{ij}^{\gamma_{ij}}\in S$ as
  $${\mathbf x}^\beta:=\prod_{l=0}^{n^2-1}x_{{\varphi^{-1}}(n^2-l)}^{\beta_{l+1}},$$
  where $\beta\in\NN^{n^2}$. For example, we write ${\mathbf x}^{(3,0,\dots,2)}$
  for $x_{nn}^2\cdot x_{21}^3$.

  Then $x^{\mathbf{\alpha}}\prec x^{\mathbf{\beta}}$ if and only if the sum of the entries
  in $\mathbf{\alpha}$ is smaller than the sum of the entries in $\mathbf{\beta}$
  or if both sums are equal, $\alpha_i>\beta_i$, where $i$ is the largest
  index in which both vectors differ.

  The essential property of $\preceq$ is the content of the following theorem. For its 
  formulation we use the notation $N(A,B)$ from Equation (\ref{nichtseiten}), see
  also Theorem \ref{Stan}.

  \begin{theorem}\label{NAB}
    For subsets $A,B\subset[n],\,\#A=\#B=k+1$ we have
    \begin{eqnarray}\lt(\det M(A,B))&=&\prod_{(i,j)\in N(A,B)}x_{ij}.\end{eqnarray}
  \end{theorem}

  For the proof we need several lemmas. In the following, let $A,B\subset [n]$
  be nonempty subsets of equal sizes.

  \begin{lemma}\label{detab} If $x_{a_{i},b_{j}}$ divides the leading monomial
    of $\det(M(A,B))$, we have that
    $$\lt(det M(A\setminus\{a_i\},B\setminus\{b_j\})\cdot
    x_{a_{i},b_{j}}=\lt(det(M(A,B)).$$
  \end{lemma}

  \begin{proof}{} We apply Laplace-expansion for the $j$-ths row
    of $M(A,B)$ and get
    \begin{eqnarray}
      \det(M(A,B))&=&\sum_{p=1}^{(k+1)}(-1)^{j+p}\cdot x_{a_p,b_j}\cdot \det(M(A\backslash
      \{a_i\},B\backslash \{b_j\}))\label{entwicklung}.
    \end{eqnarray}
    We can conclude
    $$\lt(\det(M(A,B))=\lt\Big(x_{a_i,b_j}\det(M(A\backslash\{a_i\},B\backslash\{b_j\})\Big),$$
    since $x_{a_i,b_j}\det(M(A\backslash\{a_i\},B\backslash\{b_j\}$ is the only
    summand in Equation (\ref{entwicklung}) which is a multiple of $x_{a_i,b_j}$.
    The assertion is a consequence of the fact that $\lt(mm')=\lt(m)\lt(m')$ for 
    all monomials $m,m'$.
  \end{proof}

  \begin{lemma} \label{nab} Let $(a_i,b_j)\in N(A,B)$. Then
    $N(A\setminus\{a_i\},B\setminus\{b_j\})\subset N(A,B).$
  \end{lemma}

  \begin{proof}{} Since $(a_i,b_j)\in N(A,B)$ we know that $i=(j+\ell)\mot 
    (k+1)$ for some $\ell\leq k$. We get that
    \begin{eqnarray*}
      a_{p+\ell}&>& b_p,\,\textrm{for all }p=1,\dots k+1-\ell\\
    \end{eqnarray*}
    and there is at least one $q$ such that  
    \begin{eqnarray*}
      a_{q+\ell-1}&\leq& b_q,\,q\in\{1,\dots k+1-\ell+1\}.
    \end{eqnarray*}
    Let $A':=A\backslash\{a_{j+\ell}\}=:\{a'_1<\dots<a'_{k}\},\,
    B':=B\backslash\{b_j\}=:\{b'_1<\dots<b'_{k}\}$.
    We will show that
    $$N(A',B')=N(A,B)\backslash\{(a_{(j+\ell)\mot (k+1)},b_j)\}.$$
    \begin{itemize}
    \item[(Case 1)] $j\leq k+1-\ell$: In this case 
      $\{(a'_{p+\ell},b'_p)\big|p=1,\dots,k+1\}=N(A,B)\backslash\{(a_{p+j},b_j)\}$.
      This implies 
      \begin{eqnarray}\label{Emese}a'_{p+\ell}>b'_p,\textrm{ for all }p=1,\dots,
        (k+1)-\ell-1.
      \end{eqnarray}
      We show that $\ell$ is either zero or minimal such that (\ref{Emese}) holds. 
      We get $N(A',B')=\{(a'_{p+\ell},b'_p)\big|p=1,\dots,(k+1)\}.$
      We distinguish the following cases:
      \begin{itemize}
        \item[(SubCase a)] If $j<q-1$ we have $b'_{q-1}=b_q\geq a_{q+\ell-1}=a'_{q-1+\ell-1}$, 
          thus $\ell$ is minimal in Equation (\ref{Emese}).
        \item[(SubCase b)] If $j=q-1$ then $b_j=b_{q-1}$ and $a_{j+\ell}=a_{q+\ell-1}$ are being
          removed, so that $b'_{q-1}=b_{q}$ and
          $a'_{q+\ell-2}=a_{q+\ell-2}$. Since $b_q\geq a_{q+\ell-1}$ we have $b_{q}\geq
          a_{q+\ell-2}$ and therefore $b'_{q-1}\geq a'_{q-1+\ell-1}$. We conclude that 
          $\ell$ is minimal in Equation (\ref{Emese}).
        \item[(SubCase c)]  If $j=q$ it follows that $b'_q=b_{q+1}$ and
          $a_{q+\ell-1}=a'_{q+\ell-1}$. Since $b_q\geq a_{q+\ell-1}$ it holds that 
          $b_{q+1}\geq a_{q+\ell-1}$ and in consequence we get $b'_q\geq a'_{q+\ell-1}$. Again we find
          that $\ell$ is minimal in Equation (\ref{Emese}).
          item[(SubCase d)] If $(k+1)-\ell\geq j>q$, then $b'_{q}=b_q\geq a_{q+\ell-1}=a'_{q+\ell-1}$, and
          thereby $\ell$ is minimal in Equation (\ref{Emese}).
      \end{itemize}
    \item[(Case 2)] $(k+1)\geq j>(k+1)-\ell$: Then row $a_{(j+\ell)\mot (k+1)}$ is being deleted. Since
      $(j+\ell)\mot (k+1) < p+\ell$ for $p=1,\dots,(k+1)$ we get 
      $b'_{p}=b_p,\,p=1,\dots,(k+1)-\ell$ and
      $a'_{p+\ell}=a_{p+\ell+1},\,i=1,\dots,(k+1)-\ell-1$. We conclude
      \begin{eqnarray*}
        a'_{p+\ell-1}=a_{p+\ell}&>&b_p=b'_p,\,\textrm{for all }p=1,\dots (k+1)-\ell-1 \\
      \end{eqnarray*}
      and there is at least one $q$ such that 
      \begin{eqnarray*}
        a'_{q+\ell-2}=a_{q+\ell-1}&\leq&b_q=b'_q,\,q\in\{1,\dots k+1-\ell+1\}.
      \end{eqnarray*}
      This implies 
      $$N(A',B')=\{(a'_{p+\ell-1},b'_p)\big|p=1,\dots,(k+1)-\ell-1\})=N(A,B)\backslash\{(a_{j+\ell},b_j)\}.$$
    \end{itemize}
  \end{proof}

  For the following lemma. let $\Delta^-:=\{x_{ij}\big|i>j\}$ and $\Delta^+:=\{x_{ij}\big| i\leq
  j\}$. We say that an indeterminate $x_{ij}$ is (weakly) to the right and (weakly) above
  $x_{kl}$, if $i\leq k$ and $j\geq l$. The rectangle spanned by $x_{ij}$
  and $x_{kl}$ is the set $R(x_{ij},x_{kl}):=\{x_{qr}:i\leq q\leq k,l\leq r\leq j\}.$

  Let $a,b$ be elements of $X$ such that $a$ is (weakly) to the right and (weakly) above
  $b$. We list two simple and important properties of our order $\preceq$:
  \begin{description}
    \item[(E1)] If $a,b\in\Delta^+$, then $R(a,b)$ is
      contained in $\Delta^+$ and the upper right corner is the maximal element of
      $R$ with respect to $\preceq$. The same is true for $\Delta^-$.
    \item[(E2)] If $x\in\Delta^-$ and $y\in\Delta^+$ and $x$ and $y$ share a
      column or row, we have $x\succeq y$.
  \end{description}

  \begin{lemma}\label{cap} There is an element $(i,j)\in N(A,B)$ such that
    $x_{ij}$ divides the leading monomial of $\det(M(A,B))$.
  \end{lemma}

  \begin{proof}{} Let $x_{a_{(i+\ell)\mot
    (k+1)},b_i}:=\min_{\preceq}\{x_{i,j}\big|(i,j)\in N(A,B)\}$.
    \begin{itemize}
      \item[(Case 1)] $a_{(i+\ell)\mot (k+1)}>b_i\dots$

        \begin{itemize}
          \item[(SubCase a)] $\ell=0$

            From the definition of $N(A,B)$ we get $i+\ell\leq (k+1)$.
            In this case we have that $x_{a_{(i+\ell)\mot (k+1)},b_i}$ lies in $\Delta^-$ and to the right
            and above of $x_{a_{k+1},b_1}$. Because of (\textbf{E1}) we know that
            $$V:=\{x_{a_p,b_m}\textrm{ such that }p\geq i\textrm{ and }m\leq i\}\subset \Delta^-.$$ It holds that
            $v\leq x_{a_{(i+\ell)\mot (k+1)},b_i}$ for all $v\in V$ (see Figure \ref{Fall1}).
            Assume $\det(M(A,B))$ has a term $t$ that is larger with respect 
            to $\preceq$ than the product of all $x_{ij}$ such that $(i,j)\in N(A,B)$.
            Assume further that $t$ does not contain any of those $x_{ij}$ as a factor. Then
            $t$ also contains a factor from $V$.

            If $i=1$ this assumption gives us a contradiction immediately, since then
            $t$ can not not have a factor from the first column of $M(A,B)$. Let $i>1$. 
            Since $t$ needs to have a
            factor in every row and in every column of $M(A,B)$. This means that
            $t$ comprises $i$ factors from  $W:=\{x_{a_p,b_m}$ such that $p<
            i+\ell=i$ and $m\leq i\}$. Because $W$ stretches over $i$ columns and $i-1$ rows,
            there is a row of $M(A,B)$ containing two factors of $t$. This yields the desired
            contradiction.

            \begin{figure}
              \begin{center}
                $\left(
                \begin{array}{cccccc}
                  w_{1,1} &\cdots & w_{1,i} & \star &\dots&\star \\
                  \vdots &  &\vdots &\vdots&&\vdots \\
                  w_{i-1,1} &\cdots & w_{i-1,i} & \star &\cdots &\star \\

                  v &\cdots & x_{a_i,b_i} &\star&\cdots &\star \\
                  \vdots & & \vdots &\vdots& &\vdots \\
                  v &\cdots & v &\star &\cdots&\star 
                \end{array}
                \right)$
              \end{center}
              \caption{$a_{(i+\ell)\mot (k+1)}>b_i$ and $\ell=0$}\label{Fall1}
            \end{figure}

          \item[(SubCase b)] $\ell>0$

            According to the definition of $N(A,B)$ it holds that $i+\ell\leq (k+1)$. Analogous to the
            first case, set $$V:=\{x_{a_p,b_m}\textrm{ such that  }p\geq i+\ell\textrm{ and }m\leq i\}.$$ Again
            because of (\textbf{E1}) we know $v\leq x_{a_{(i+\ell)\mot (k+1)},b_i}$ for all $v\in V$.
            Since $\ell>0$ there is a $j\leq (k+1)-\ell+1$ such that  $a_{j+\ell-1}\leq b_j$ and so we know
            $x_{a_{j+\ell-1},b_j}\in\Delta^+$ (see Figure \ref{Fall2}).

            We consider the following subcases:
            \begin{itemize}
              \item[$j \leq i$:] 

                Here $x_{a_i,b_i}$ lies to the right and below of
                $x_{a_{j+\ell-1},b_j}$ and is therefore an element of $\Delta^+$ as well. This
                together with (\textbf{E1}) is the reason why the rectangle
                $$V':=\{x_{a_p,b_m}\textrm{ such that }p\leq j+\ell-1\textrm{ and }m\leq i\}$$ is
                contained in $\Delta^+$. By (\textbf{E2}) we know that $x_{a_1,b_i}\prec
                x_{a_{(i+\ell)\mot (k+1)},b_i}$. Again (\textbf{E1}) tells us that $v'\prec
                x_{a_1,b_i}\prec x_{a_{(i+\ell)\mot (k+1)},b_i}$ for all $v'\in V'$.

                Assume $\det(M(A,B))$ had a monomial $t$ that was larger w.r.t. $\preceq$
                than the product of all $x_{i,j}$ such that  $(i,j)\in N(A,B)$ and none of those
                $x_{ij}$ was a factor of $t$. Then $t$ would not contain a factor from $V\cup V'$.
                The set $V\cup V'$ stretches over $i-j+1$ columns of  $M(A,B)$
                inside of which all except for $i+\ell-(j+\ell-1)-1=i-j$ rows are covered.
                This means, that $t$ had two factors in a row or, if $i=j$, no factor in a
                certain column, which is a contradiction.

                \begin{figure}
                  \begin{center}
                    $\left(
                    \begin{array}{cccccc}
                      \star & v' & v' & \star &\star &\star \\
                      \star & v' & v' &\star &\star &\star \\
                      \star & x_{a_{j+\ell-1},b_j} & v' & \star &\star &\star \\

                      \star &\star & \star &\star&\star &\star \\
                      v &v & x_{a_{i+\ell},b_i} &\star&\star&\star  \\
                      v &v & v &\star &\star&\star \end{array}\right)$ $\left(\begin{array}{cccccc}
                      \star &\star & \star &\star&\star &\star \\
                      \star &\star & \star &\star&\star &\star \\
                      \star &\star & \star &\star&\star &\star \\
                      v & x_{a_{i+\ell},b_i}&\star &\star&v'&v'  \\
                      v &v &\star &\star&v'&v'  \\
                      v &v & \star &\star &x_{a_{j+\ell-1},b_j}&v' 
                    \end{array}
                    \right)$
                  \end{center}
                  \caption{$a_{(i+\ell)\mot (k+1)}>b_i$ and $(k+1)\geq \ell>0$}\label{Fall2}
                \end{figure}

              \item[$j >i$:] $j>i$

                Here we have that $x_{a_{i+\ell},b_m}$ lies to the right and above
                $x_{a_{j+\ell-1},b_j}$ and is therefore contained in $\Delta^+$. Property
                (\textbf{E1}) tells us that the rectangle
                $$V':=\{x_{a_p,b_m}\textrm{ such that }i+\ell\leq p\leq j+\ell-1\textrm{ and }m\geq
                j\}$$ is a subset of $\Delta^+$. Because of (\textbf{E2}) we know
                $x_{a_{i+\ell},b_{k+1}}\prec x_{a_{i+\ell},b_i}$ and again (\textbf{E1}) tells us that
                $v'\preceq x_{a_{(i+\ell)\mot (k+1)},b_i}$ for all $v'\in V'$.

                Again we assume that $\det(M(A,B))$ had a term $t$ that was larger than the
                product of all $x_{ij}$ such that $(i,j)\in N(A,B)$ and none of those $x_{ij}$
                was a factor of $t$. Again $t$ could not contain a factor from $V\cup V'$.
                Since $V\cup V'$ comprises $j+\ell-1-(i+\ell)+1=j-i$ rows of $M(A,B)$
                inside of which all except of $j-i-1$ columns are covered.
                That means that $t$ had two factors in a column, which is again a
                contradiction.
            \end{itemize}
        \end{itemize}

      \item[(Case 2)] $a_{(i+\ell)\mot (k+1)}>b_i$ and $(k+1)\geq \ell>0$.

        According to the definition of $N(A,B)$ we have $i+\ell>(k+1)$ in this case and
        there is a $j\leq (k+1)-\ell+1$ such that
        $a_{j+\ell-1}<b_j$. Thus we know $x_{a_{j+\ell-1},b_j},x_{a_{i+\ell},b_i}\in
        \Delta^+$ and
        $x_{a_{i+\ell},b_i}$ is to the right and above $x_{a_{j+\ell-1},b_j}$. For the
        enclosed rectangle $$V:=\{x_{a_p,b_m}\textrm{ such that } j+\ell-1\geq p\geq
        (i+\ell)\mot (k+1)\textrm{ and } j\leq m\leq i\},$$ we know that $v\leq
        x_{a_{(i+\ell)\mot (k+1)},b_i}$ for all $v\in V$ (see Figure \ref{Fall3}).

        Again we assume $\det(M(A,B))$ had a $t$ that was larger with respect to our
        term-order than the product of all $x_{ij}$ such that $(i,j)\in N(A,B)$ and that
        did not contain one of the $x_{ij}$ as a factor. In consequence $t$ could
        not have a factor from $V$.
        For $i>j$ the set $V$ stretches over  $j+\ell-1-(i+\ell)\mot
        (k+1)\,+1=j+\ell-(i+\ell-(k+1))=(k+1)-i+j$ rows of $M(A,B)$
        and covers $i-j+1$ columns of it, leaving only $(k+1)-(i-j+1)=(k+1)-i+j-1$ uncovered.
        For $i=j$ we know that $V$ consists of the full $i$th column of $M(A,B)$.
        That means that $t$ contains either two or no factor in a column, a
        contradiction!
    \end{itemize}
  \end{proof}

  \begin{figure}
    \begin{center}
      $\left(
      \begin{array}{cccccc}
        \star &\star & \star &\star&\star &\star \\
        \star &v & v &v&x_{a_{i+\ell},b_i} &\star \\
        \star &v & v &v&v&\star \\
        \star &x_{a_{j+\ell-1},b_j} & v &v&v &\star \\
        \star &\star & \star &\star&\star &\star \\
        \star &\star & \star &\star&\star &\star 
      \end{array}
      \right)$
    \end{center}
    \caption{$a_{(i+\ell)\mot (k+1)}>b_i$ and $(k+1)\geq \ell>0$}\label{Fall3}
  \end{figure}

  \begin{proof}[Proof of Theorem \ref{NAB}]
    For $k=0$ the statement is clear. For $k>0$ we proceed via an induction. The leading monomial of
    $\det M(A,B)$ contains according to Lemma \ref{cap} a factor $x_{a_i,b_j}\in
    N(A,B)$. A consequence of Lemma \ref{detab} is the fact, that
    $$\lt(\det M(A,B))=x_{a_i,b_j}\cdot\lt(\det M(A',B')),$$
    where $A'=A\backslash\{a_i\},\,B'=B\backslash\{b_j\}$. The induction's
    assumption together with Lemma \ref{nab} yields the theorem.
  \end{proof}

\end{document}